\documentclass[12pt]{amsart}
\usepackage{graphicx}
\usepackage{latexsym,amscd,amssymb,epsfig} 

\DeclareMathOperator{\sd}{{\rm sd}}
\DeclareMathOperator{\des}{{\rm des}}
\DeclareMathOperator{\Hilb}{{\rm Hilb}}

\def\RR{{\mathbb R}}

\def\fvec{{\mathfrak f}}
\def\fpol{{\mathfrak f}}
\def\hvec{{\mathfrak h}}
\def\hpol{{\mathfrak h}}

\def\fA{{\mathfrak A}}

\def\xx{{\bf {\underline{x}}}}
\def\fmat{{\mathfrak F}}
\def\hmat{{\mathfrak H}}

\theoremstyle{plain}
  \newtheorem{theorem}{Theorem}[section]
  
  \newtheorem{lemma}[theorem]{Lemma}
  \newtheorem{corollary}[theorem]{Corollary}
  
\theoremstyle{definition}
  
  \newtheorem{example}[theorem]{Example}
  
  \newtheorem{question}[theorem]{Question}
  \newtheorem{remark}[theorem]{Remark}
  \newtheorem{problem}[theorem]{Problem}

\numberwithin{equation}{section}

\begin{document}

\title[$f$-Vectors of Barycentric Subdivisions]
{$f$-Vectors of Barycentric Subdivisions}

\author{Francesco Brenti}
\address{Dipartimento di Matematica \\
Universita' di Roma "Tor Vergata" \\
Via della Ricerca Scientifica, 1 \\
00133, Roma, Italy}

\email{brenti@mat.uniroma2.it}

\author{Volkmar Welker}
\address{Fachbereich Mathematik und Informatik\\
Philipps-Universit\"at Marburg\\
35032 Marburg, Germany}
\email{welker@mathematik.uni-marburg.de}

\thanks{Both authors partially supported by EU
      Research Training Network 
``Algebraic Combinatorics in Europe'', grant HPRN-CT-2001-00272
and the program on ``Algebraic Combinatorics'' at the Mittag-Leffler Institut in
Spring 2005}



\begin{abstract}
For a simplicial complex or more generally Boolean cell complex $\Delta$ 
we study the behavior of the $f$- and
$h$-vector under barycentric subdivision. We show that if $\Delta$ has
a non-negative $h$-vector then the $h$-polynomial of its barycentric 
subdivision has only simple and real zeros. As a consequence this implies 
a strong version of the Charney-Davis conjecture for spheres that are the 
subdivision of a Boolean cell complex.  

For a general $(d-1)$-dimensional
simplicial complex $\Delta$ the $h$-polynomial of its $n$-th iterated
subdivision shows convergent behavior. More precisely, we show that
among the zeros of this $h$-polynomial there is one converging
to infinity and the other $d-1$ converge to a set of $d-1$ real numbers
which only depends on $d$. 
\end{abstract}

\maketitle

\section{Introduction}
\label{introduction}

This work is concerned with the effect of barycentric subdivision on the 
enumerative structure of a simplicial complex. More precisely, we study the
behavior of the $f$- and $h$-vector of a simplicial complex, or more
generally, Boolean cell complex, under barycentric subdivisions.
In our main result we show that if the $h$-vector of a simplicial complex 
or Boolean cell complex is non-negative
then the $h$-polynomial (i.e., the generating polynomial of the $h$-vector)
of its barycentric subdivision has only  real simple
zeros. 
Moreover, if one applies barycentric subdivision iteratively then there
is a limiting behavior of the zeros of the $h$-polynomial, with one zero going to infinity
and the other zero converging. The limit values of the zeros only depend on the 
dimension of the complex.

In this introductory section we first review the basic concepts around
$f$- and $h$-vectors and give the basic definitions, then we outline the
structure of this work.

Our results hold in the generality of Boolean cell complexes.
Recall that if $\Delta$ is a CW-complex then for two open cells $A$ and $A'$ in
$\Delta$ we define $A \leq_\Delta A'$ if $A$ is contained in the closure of $A'$ in
$\Delta$. If $\Delta$ is a regular CW-complex then this partial order on
its cells already encodes the topology of $\Delta$ up to homeomorphism.
In this text we consider the empty cell as a cell of any CW-complex
$\Delta$ which then serves as the least element of the partial order
on $\Delta$. 
Thus we can identify $\Delta$ with the partial order on its cells. 
In particular, we write $A \in \Delta$ if we want to express that 
$A$ is a cell of $\Delta$. 
We call a regular CW-complex $\Delta$ a Boolean cell complex if 
for each $A \in \Delta$ the lower interval $[\emptyset,A] :=
\{ B \in \Delta ~|~\emptyset \leq_\Delta B \leq_\Delta A\}$ is a Boolean lattice -- i.e.
the lattice of subsets of a set. The prime example of a Boolean cell
complex is a simplicial complex. 
Most of the time we consider a
simplicial complex as an abstract simplicial complex and therefore
identify its open cells with subsets of the ground set of $\Delta$. 
Adopting notions from simplicial complexes, we call a
cell $A \in \Delta$ a face of $\Delta$.  
By $\dim (A)$ we denote the dimension of the cell $A$.
The dimension $\dim(\Delta)$ of $\Delta$ is the maximal dimension of one
of its faces. For $-1 \leq i \leq d-1$ ($\stackrel{\rm def}{=}
\dim(\Delta)$) we write $f_i^\Delta$ for
the number of $i$-dimensional faces of $\Delta$.
The vector $\fvec^\Delta = (f_{-1}^\Delta, \ldots, f_{d-1}^\Delta)$ 
is called the $f$-vector of $\Delta$. 
Often it will be convenient to encode the $f$-vector of a Boolean cell complex
as a polynomial. By $\fpol^\Delta(t) = \sum_{i=0}^d
f_{i-1}^\Delta t^{d-i}$ we 
denote the $f$-polynomial of the $(d-1)$-dimensional Boolean cell complex
$\Delta$. Information equivalent to the one carried by the $f$-vector is
encoded in the $h$-vector. It is defined as $\hvec^\Delta = (h_0^\Delta,
\ldots, h_d^\Delta)$ where $\hpol^\Delta (t) = \sum_{i=0}^d
h_i^\Delta t^{d-i}$ is the expansion of $\fpol^\Delta(t-1)$ in terms of $t$-powers.

The barycentric subdivision of a Boolean cell complex $\Delta$ is the simplicial
complex $\sd(\Delta)$, that as an abstract simplicial complex
is defined on the ground set $\Delta \setminus \{ \emptyset\}$ -- the set of
non-empty faces of $\Delta$ --  with 
$i$-faces the strictly increasing flags $A_0 <_\Delta \cdots <_\Delta A_i$ of faces 
in $\Delta \setminus \{ \emptyset \}$. It is well known that 
$\Delta$ and $\sd(\Delta)$ are homeomorphic. Thus 
$\Delta$ and $\sd(\Delta)$ define cellulations and  triangulations
of the same space. 

Throughout the paper we will use $[n]$ to denote the set $\{1, \ldots, n\}$ for a natural
number $n \geq 0$.

\section{The $h$-Vector of a Barycentric Subdivision}
\label{h-vector-section}

For the sake of completeness we sketch the proof of the following well known result.

\begin{lemma} \label{f-subdivision} Let $\Delta$ be a $(d-1)$-dimensional Boolean cell complex then
$$f_j^{\sd(\Delta)} = \sum_{i = 0}^d f^{\Delta}_{i-1} \, (j+1)! \, S(i,j+1),$$
where $S(j,i)$ is the Stirling number of the second kind.
\end{lemma}
\begin{proof}  
By definition a $j$-face of $\sd(\Delta)$ is a flag $A_0 <_\Delta \cdots <_\Delta A_j$ 
of faces in $\Delta \setminus \{\emptyset\}$. Let us fix $j$ and $A_j$. Since $\Delta$ is
a Boolean cell complex, the interval $[\emptyset,A_j]$ is the Boolean lattice of subsets
of a $(\dim(A_j)+1)$-element set. Thus we may identify each of the cells $A_0, \ldots, A_j$ with
a subset of this set. Using this identification, the map sending $A_0 <_\Delta \cdots <_\Delta  A_j$ 
to $(A_0, A_1 \setminus A_0 , \cdots , A_j \setminus
A_{j-1})$ is a bijection between $j$-faces of $\sd(\Delta)$ with $A_j$ as its largest
element with respect to $\leq_\Delta$ and ordered partitions of $A_j$ into $j+1$ non-empty 
blocks. The latter are enumerated by
$(j+1)! S(i,j+1)$, where $i = \dim(A_j)+1$.
Now summing over all $A_j$ of fixed dimension $i-1$ we get $f^{\Delta}_{i-1} (j+1)! S(i,j+1)$. Summing
over all $i$ then yields the desired formula.
\end{proof}

Less obvious is the representation of the $h$-vector of a barycentric 
subdivision in terms of the $h$-vector of the original complex.
In the formulation of the proposition we write 
$$
D(\sigma)=\{ i \in [d-1] : \, \sigma(i) > \sigma(i+1) \},
$$
$\des(\sigma)=\#D(\sigma)$ for the number of 
descents of the permutation $\sigma$, and $S_d$ for the symmetric group on $[d]$.
For $1 \leq d$, $1 \leq j \leq d$ and $0 \leq i \leq d-1$ ,
we denote by $A(d,i,j)$ the number of permutations $\sigma \in S_d$ such that 
$\sigma(1) = j$ and $\des(\sigma) = i$. We define $A(d,i,j)$ for all $d \geq 1$ and all 
integers $i$ and $j$. In particular, $A(d,i,j) = 0$ if $i \leq -1$ or $i \geq d$.
A formula similar to the one in the following theorem appears in a slightly different
context in \cite[Chapter 3, Ex. 71]{st1} and \cite[Theorem 8.3]{st2}.
We were not able to deduce the following result from these formulas directly.

\begin{theorem}
\label{2.2}
\label{h-subdivision}
Let $\Delta$ be a $(d-1)$-dimensional Boolean cell complex. Then
$$h_j^{\sd(\Delta)} = \sum _{r=0}^{d} A(d+1,j,r+1) \, h_r^\Delta$$
for $0 \leq j \leq d$.
\end{theorem}
\begin{proof} 
For all $0 \leq j \leq d$ we have that
\begin{eqnarray*}
h_{j}^{sd(\Delta )} & = &  
\sum _{i=0}^{j} \left( ^{^{\scriptstyle d-i}}_{_{\scriptstyle
j-i}}\right) (-1)^{j-i} f_{i-1}^{sd (\Delta )} \\
& =& \sum _{i=0}^{j} \left( ^{^{\scriptstyle d-i}}_{_{\scriptstyle
j-i}}\right)  (-1)^{j-i} \, \sum _{k=0}^{d} f_{k-1}^{\Delta } \, S(k,i) \, i! \\
& = & \sum _{i=0}^{j} \sum _{k=0}^{d} \left( ^{^{\scriptstyle d-i}}_{_{\scriptstyle
j-i}}\right) (-1)^{j-i} \, S(k,i) \, i! \, \sum _{r=0}^{k}
\left( ^{^{\scriptstyle d-r}}_{_{\scriptstyle
d-k}}\right) h_{r}^{\Delta } \\
& = & \sum _{r=0}^{d} \left( \sum _{k=0}^{d} \sum_{i=0}^{j} (-1)^{j-i}
\left( ^{^{\scriptstyle d-i}}_{_{\scriptstyle
j-i}}\right) \left( ^{^{\scriptstyle d-r}}_{_{\scriptstyle
d-k}}\right) \, S(k,i) \, i! \right) \, h_{r}^{\Delta }.
\end{eqnarray*}

Now notice that if $S = \{ s_{1}, \ldots , s_{k} \} \subseteq [d]$ and $s_{1}<\ldots <s_{k}$
then
\[ \#\Big\{ \sigma \in S_{d+1}: \; {\small \begin{array}{c} D(\sigma ) \subseteq S, \\
\sigma (d+1) =d+1-r \end{array}} \Big\} =
\left( ^{^{\scriptstyle \hspace{1.2cm} s_{k}}}_{_{\scriptstyle
s_{k}-s_{k-1},\ldots ,s_{2} -s_{1},s_{1}}}\right)
\left( ^{^{\scriptstyle \hspace{0.1cm} d-r}}_{_{\scriptstyle d-s_{k}}}\right). \]

So,
\begin{eqnarray*} 
&   &  \sum_{\{S \subseteq [d], \; \#S=k \} }
\#\Big\{ \sigma \in S_{d+1}: \; {\small \begin{array}{c} D(\sigma ) \subseteq S, \\ \sigma (d+1)
=d+1- r\end{array}} \Big\} \\
& = & \sum _{1 \leq s_{1} < \ldots < s_{k} \leq d} \left( ^{^{\scriptstyle \hspace{1.2cm}
 s_{k}}}_{_{\scriptstyle
s_{k}-s_{k-1},\ldots ,s_{2}-s_{1},s_{1}}}\right)
\left( ^{^{\scriptstyle \hspace{0.1cm} d-r}}_{_{\scriptstyle d-s_{k}}}\right) \\
& = & \sum _{j=k}^{d} \left( ^{^{\scriptstyle d-r}}_{_{\scriptstyle d-j}}\right)  \,
\sum_{1\leq s_{1}< \ldots <s_{k-1} \leq j-1}
 \left( ^{^{\scriptstyle \hspace{1.2cm} j}}_{_{\scriptstyle
j-s_{k-1},\ldots ,s_{2}-s_{1},s_{1}}}\right) =
 \sum_{j=k}^{d} \left( ^{^{\scriptstyle d-r}}_{_{\scriptstyle d-j}}\right) S(j,k) \, k!
\end{eqnarray*} 

Therefore, 
\[ \sum_{k=0}^{d} \sum_{i=0}^{j} (-1)^{j-i} 
\left( ^{^{\scriptstyle d-i}}_{_{\scriptstyle
j-i}}\right) \left( ^{^{\scriptstyle d-r}}_{_{\scriptstyle
d-k}}\right) \, S(k,i) \, i! \]

\[ = \sum_{i=0}^{j}(-1)^{j-i} \left( ^{^{\scriptstyle d-i}}_{_{\scriptstyle
j-i}}\right)   \sum_{\{S\subseteq [d], \, \#S=i \}} \#\Big\{ \sigma \in S_{d+1}:
{\small \begin{array}{c} D(\sigma ) \subseteq  S, \\ \sigma (d+1)=d+1-r \end{array}} \Big\}  \]

\[ = \sum_{\{ S \subseteq [d], \, \#S \leq j\} }(-1)^{j-\#S} \,
\left( ^{^{\scriptstyle d-\#S}}_{_{\scriptstyle
j-\#S}}\right) \, \#\Big\{ \sigma \in S_{d+1}:  {\small \begin{array}{c} D(\sigma )\subseteq S, \\
\sigma (d+1)=d+1-r \end{array}} \Big\}  \]

\[ = \sum_{\{ S \subseteq [d], \, \#S \leq j\} }(-1)^{j-\#S} \,
\left( ^{^{\scriptstyle d-\#S}}_{_{\scriptstyle j-\#S}}\right) 
\sum_{T \subseteq S} \#\Big\{ \sigma \in S_{d+1}: {\small \begin{array}{c}D(\sigma )=T, \\ \sigma (d+1)
=d+1-r \end{array}} \Big\} \]

\[ = \sum_{\{ T \subseteq [d], \; \#T \leq j \} } \#\Big\{ \sigma \in S_{d+1}:
{\small \begin{array}{c} D(\sigma )
=T, \\ \sigma (d+1)=d+1-r \end{array}} \Big\} \, \sum _{\{ S \supseteq T, \; \#S \leq j\}}(-1)^{j-\#S}
\, \left( ^{^{\scriptstyle d-\#S}}_{_{\scriptstyle j-\#S}}\right)   \]

\[ = \sum_{\{ T \subseteq [d], \; \#T \leq j\} } \# \Big\{ \sigma \in S_{d+1}: {\small \begin{array}{c} D(\sigma )
=T, \\ \sigma (d+1)=d+1-r \end{array}} \Big\} \, \sum_{i=\#T}^{j} (-1)^{j-i}\left( ^{^{\scriptstyle d-i}}_{_{\scriptstyle
j-i}}\right)   \left( ^{^{\scriptstyle d-\#T}}_{_{\scriptstyle
i-\#T}}\right). \]
But
\begin{eqnarray*}
\sum_{i=\#T}^{j}(-1)^{j-i}\left( ^{^{\scriptstyle d-i}}_{_{\scriptstyle j-i}} \right)
\left( ^{^{\scriptstyle d-\#T}}_{_{\scriptstyle i-\#T}} \right) & = &
\left( ^{^{\scriptstyle d-\#T}}_{_{\scriptstyle j-\#T}} \right) \,
\sum_{i=\#T}^{j}(-1)^{j-i} \, \left( ^{^{\scriptstyle j-\#T}}_{_{\scriptstyle i-\#T}} \right) \\
& = & \delta_{j,\#T}
\end{eqnarray*}
Hence
\[ \sum_{k=0}^{d} \, \sum_{i=0}^{j}(-1)^{j-i} \,
\left( ^{^{\scriptstyle d-i}}_{_{\scriptstyle j-i}} \right) \,
\left( ^{^{\scriptstyle d-r}}_{_{\scriptstyle d-k}} \right) \, S(k,i) i! \]

\[ = \sum_{\{ T \subseteq [d], \; \#T=j \} } \#\Big\{ \sigma \in S_{d+1}: {\small \begin{array}{c} D(\sigma )
=T, \\ \sigma (d+1)=d+1-r \end{array}} \Big\} \]

\[ = \#\Big\{ \sigma \in S_{d+1}: {\small \begin{array}{c} \des(\sigma)=j,
 \\ \sigma(d+1)=d+1-r \end{array}} \Big\} , \]
and the result follows.
\end{proof}

\begin{corollary} \label{monotone}
Let $\Delta$ be a $(d-1)$-dimensional Boolean cell complex such that
$h_i^\Delta \geq 0$ for $1 \leq i \leq d$. Then for $0 \leq i \leq d$
$$h_i^{\sd(\Delta)} \geq h_i^\Delta.$$
\end{corollary}
\begin{proof} We have from Theorem \ref{2.2} and our hypotheses that
$h_{j}^{sd(\Delta)} \geq A(d+1, j,j+1)\, h_{j}^{\Delta}$ so the result
follows since $A(d+1,j,j+1)\geq1$ for all $0\leq j \leq d $.
\end{proof}

\begin{example} \label{example-nonincreasing}
It is not true that $h_i^{\sd(\Delta)} \geq h_i^\Delta$ for any simplicial complex $\Delta$ 
and any $i$. For example let $\Delta$ be the disjoint union of a $2$-simplex and a triangle
with $f$-vector $\fvec^\Delta = (1,6,6,1)$ and $h$-vector $\hvec^\Delta = (1,3,-3,0)$.
Then $\fvec^{\sd(\Delta)} = (1,13,18,6)$ and $\hvec^{\sd(\Delta)} = (1,10,-5,0)$.
\end{example}

Before we proceed with the next corollary we summarize some basic facts about the numbers
$A(d,i,j)$.

\begin{lemma} \label{basiclemma} ~

\begin{itemize}
\item[(i)] 
$$A(d,i,j) = \displaystyle{\sum_{l=1}^{j-1} A(d-1,i-1,l) + \sum_{l=j+1}^{d}} A(d-1,i,l-1),$$
for $2 \leq d$, $1 \leq j \leq d$, and $0 \leq i \leq d-1$.

\item[(ii)] 
$$A(d,i,j) = A(d,d-1-i,d+1-j)$$ for $1 \leq d$, 
$1 \leq j \leq d$ and $0 \leq i \leq d-1$.
\end{itemize}
\end{lemma}
\begin{proof} 
For (i) the recursion formula just reflects the enumeration of the permutations 
in $S_d$ with $\sigma(1) = j$ and $\des(\sigma) = i$ by $l = \sigma(2)$.
For (ii) the assertion reflects the fact that the map $\sigma (1) \cdots \sigma (d)
\leftrightarrow d+1-\sigma (1) \cdots d+1-\sigma (d)$ is a bijection between the
sets enumerated by the the two numbers.
\end{proof}

We call the $h$-vector $\hvec^\Delta$ (resp. $h$-polynomial $h^{\Delta}(t)$)
of a Boolean cell complex
$\Delta$ reciprocal if $h_i^\Delta = h_{d-i}^\Delta$ for $0 \leq i \leq d$ 
(resp. $\hpol^\Delta(t) = t^d \hpol^\Delta(\frac{1}{t})$). Clearly, $\hvec^\Delta$
is reciprocal if and only if $\hpol^\Delta(t)$ is. Either condition is
equivalent to the $h$-vector satisfying the Dehn-Sommerville relations (see \cite[Theorem 5.4.2]{BH}).
 
\begin{corollary} \label{reciprocal}
Let $\Delta$ be a $(d-1)$-dimensional Boolean cell complex with reciprocal $h$-vector
then $\sd(\Delta)$ also has reciprocal $h$-vector.
\end{corollary}
\begin{proof}
By Lemma \ref{basiclemma} and Theorem \ref{h-subdivision} we get
\begin{eqnarray*}
h_j^{\sd(\Delta)} & = & \sum_{i=0}^d A(d+1,j,i+1) h_i^\Delta.\\
                  & = & \sum_{i=0}^d A(d+1,d-j,d-i+1) h_{d-i}^\Delta \\
                  & = & \sum_{i=0}^d A(d+1,d-j,i+1) h_i^\Delta \\
                  & = & h_{d-j}^{\sd(\Delta)} \\
\end{eqnarray*}
\end{proof}

\section{Main Result}
\label{main-result}

For the formulation of the following result we recall that a sequence $(a_0,\ldots,
a_d)$ of real numbers is called log-concave if $a_i^2 \geq a_{i-1}a_{i+1}$ for $1 \leq i \leq d-1$.
We say that $(a_0,\ldots, a_d)$ has internal zeros, if there are $0 \leq i < j \leq d$ such that
$a_i, a_j \neq 0$ and $a_{i+1} = a_{i+2} = \cdots = a_{j-1} = 0$.
It is a well known and easy to prove fact that if $(a_0,\ldots,
a_d)$ is a log-concave
sequence of non-negative real numbers and has no internal zeros then it is unimodal;
that is there is a $0 \leq j \leq d$ such that $a_0 \leq \cdots \leq a_j \geq \cdots \geq a_d$.
By another well known fact (see, e.g., \cite[\S 7.1]{})
it follows that if the polynomial $a_0 + a_1t + \cdots + a_dt^d$
has only real zeros then $(a_0,\ldots, a_d)$ is log-concave.

\begin{theorem}
\label{realroots}
Let $\Delta$ be a $(d-1)$-dimensional Boolean cell complex such that $h_i^\Delta 
\geq 0$ for $0 \leq i \leq d$. Then
$$\hpol^{\sd(\Delta)}(t) = \sum_{i=0}^d h_i^{\sd(\Delta)} t^{d-i}$$ has only simple and
real zeros. In particular, 
$\hvec^{\sd(\Delta)} = (h_0^{\sd(\Delta)}, \ldots, h_d^{\sd(\Delta)})$ is a log-concave 
sequence. If in addition the sequence $\hvec^{\sd(\Delta)}$ has no internal zeros  then
it is unimodal.
\end{theorem}

Before we can prove Theorem \ref{realroots} we need a few definitions and a 
 lemma about permutation statistics of Coxeter groups of types $A_d$ and $B_d$.
The Coxeter group of type $A_d$ is the symmetric group $S_{d+1}$ and we have defined the respective
statistics already in Section \ref{h-vector-section}. The Coxeter group of type $B_d$ is the
group of all bijections $\sigma$ of $ \{ \pm 1, \ldots, \pm d\}$ in itself such that 
$\sigma(-i) = - \sigma(i)$ for $1 \leq i \leq d$. An element of $B_d$ is called a
signed permutation. For $\sigma \in B_d$ one conveniently sets $\sigma(0) := 0$. 
Then $i \in [0,d-1]$ is a descent of $\sigma \in B_d$ if $\sigma(i) > \sigma(i+1)$. 
In particular, $0$ is a descent if and only if $\sigma(1) < 0$. 
By $\des_B(\sigma)$ we denote the number of
descents of the signed permutation $\sigma$. By $N(\sigma)$ we denote the number of 
$i \in [d]$ such that $\sigma(i) < 0$.

\begin{lemma} \label{bijection} 
Let $d \geq 1$, $j \in [d+1]$ and $0 \leq i \leq d$. There is a bijection between pairs $(S,\sigma)$ where $S$ is a 
$(j-1)$-subsets of $[d]$ and $\sigma \in S_{d+1}$ has $i$ descents and $\sigma (1) = j$ and
the set of signed permutations $\tau \in B_d$
with $i$ descents and $N(\tau) = j-1$.
\end{lemma}
\begin{proof} 
Let $\fA(d+1,i,j)$ be the set of pairs $(S,\sigma)$ of $(j-1)$-subsets
$S$ of $[d]$ and $\sigma \in S_{d+1}$ such that
$\sigma$ has $i$ descents and $\sigma(1) = j$. Then for $(S,\sigma) \in \fA(d+1,i,j)$
we define $\tau_{(S,\sigma)} \in B_d$ as follows. Let $S = \{ s_1 > \ldots > s_{j-1} \}$ and
$[d] \setminus S = \{ s_j < \ldots < s_d\}$. Then define $\tau_{(S,\sigma)} \in B_d$ by

$$\tau_{(S,\sigma)}(i) = \left\{ \begin{array}{ccc} -s_{\sigma(i+1)} & \mbox{~if~} & \sigma(i+1) < j \\
                                               s_{\sigma(i+1)-1} & \mbox{~if~} & \sigma(i+1) > j
\end{array} \right. .$$

We show that for $0 \leq i \leq d$:

\begin{eqnarray} \tau_{(S,\sigma)} (i) > \tau_{(S,\sigma)}(i+1) \Leftrightarrow \sigma(i+1) >
 \sigma(i+2) .
\end{eqnarray}

This is clear for $i=0$, so assume $1 \leq i \leq d$.
Let us distinguish the four possible cases:

\begin{itemize}
\item[$\triangleright$] $\sigma(i+1), \sigma(i+2) < j$.
In this case $\tau_{(S,\sigma)} (i) > \tau_{(S,\sigma)}(i+1) \Leftrightarrow s_{\sigma(i+1)}
< s_{\sigma(i+2)}$,
which is equivalent to $\sigma(i+1) > \sigma(i+2)$.
\item[$\triangleright$] $\sigma(i+1), \sigma(i+2) > j$.
In this case $\tau_{(S,\sigma)} (i) > \tau_{(S,\sigma)}(i+1) \Leftrightarrow s_{\sigma(i+1)-1} >
 s_{\sigma(i+2)-1}$, which is equivalent to $\sigma(i+1) > \sigma(i+2)$.
\item[$\triangleright$] $\sigma(i+1) < j,  \sigma(i+2) > j$.
In this case clearly $\sigma(i+1) < \sigma(i+2)$ and $\tau_{(S,\sigma)} (i) < 0 <  \tau_{(S,\sigma)}(i+1)$.
\item[$\triangleright$] $\sigma(i+1) > j,  \sigma(i+2) < j$.
In this case clearly $\sigma(i+1) > \sigma(i+2)$ and $\tau_{(S,\sigma)} (i) > 0 > \tau_{(S,\sigma)}(i+1)$.
\end{itemize}

Thus $\des_B(\tau_{(S,\sigma)}) = \des(\sigma)$. Clearly $N(\tau_{(S,\sigma)}) = j-1$. One easily checks 
that the map $(S,\sigma) \mapsto \tau_{(S,\sigma)}$ is injective and surjective from $\fA(d+1,i,j)$ to
the set of signed permutations in $B_d$ with $i$ descents and $N(\tau) = j-1$.
\end{proof}

\begin{proof}[Proof of Theorem \ref{realroots}]  
By Theorem \ref{h-subdivision} the following identity holds.

\begin{eqnarray}
\label{first} \frac{\displaystyle{\sum_{i=0}^d} h_i^{\sd(\Delta)}t^{i}}{(1-t)^{d+1}} & = & 
{\sum_{j=0}^d} h_j^\Delta \cdot \frac{\displaystyle{\sum_{i=0}^d} A(d+1,i,j+1) t^{i}}{(1-t)^{d+1}} 
\end{eqnarray}

If $B(d,i,j)$ is the number of signed permutations in $B_d$ such that $\des_B(\sigma) = i$ and $N(\sigma) = j$
then it follows from \cite[Corollary 3.9, eq.(25)]{Brenti-toric} that:

\begin{eqnarray}
\label{second} \frac{\displaystyle{\sum_{i=0}^d} B(d,i,j) t^{i}}{(1-t)^{d+1}} & = &
\sum_{r \geq 0} {d \choose j} r^j(r+1)^{d-j} t^r 
\end{eqnarray}

Using Lemma \ref{bijection} it follows from  (\ref{second}) that

\begin{eqnarray}
\label{third} \frac{\displaystyle{\sum_{i=0}^d} A(d+1,i,j+1) t^{i}}{(1-t)^{d+1}} & = &
\sum_{r \geq 0} r^j(r+1)^{d-j} t^r.
\end{eqnarray} 

Therefore, by (\ref{first})

\begin{eqnarray}
\label{fourth}
\frac{\displaystyle{\sum_{i=0}^d} h_i^{\sd(\Delta)}t^{i}}{(1-t)^{d+1}} & = & 
\sum_{r \geq 0}( {\sum_{j=0}^d} h_j^\Delta r^j(r+1)^{d-j} ) \, t^r.
\end{eqnarray}

By  \cite[Theorem 4.2]{Braenden} it follows that if we write

\begin{eqnarray}
\label{fifth}
 \sum_{j=0}^d h_j^\Delta r^j(r+1)^{d-j}  = \sum_{j=0}^d k_j^\Delta {r \choose j}
\end{eqnarray}
then the polynomial $\displaystyle{\sum_{j=0}^d k_j^\Delta r^j}$ has nonnegative coefficients 
and only simple and real zeros.  Using the Binomial Theorem we obtain, by (\ref{fourth})
and (\ref{fifth})

\begin{eqnarray*}
\frac{\displaystyle{\sum_{i=0}^d} h_i^{\sd(\Delta)}t^{i}}{(1-t)^{d+1}}
& = & \displaystyle{\sum_{r \geq 0}}( \displaystyle{\sum_{j=0}^d} k_j^\Delta {r \choose j} )
\, t^r   \\
& = & \displaystyle{\sum_{j=0}^d}  k_j^\Delta  \sum_{r \geq 0} {r \choose j} \, t^r    \\
& = & \displaystyle{\sum_{j=0}^d}  k_j^\Delta  \frac{t^j}{(1-t)^{j+1}}    \\
& = & \displaystyle{\frac{\displaystyle{\sum_{j=0}^d  k_j^\Delta t^j (1-t)^{d-j} }}{(1-t)^{d+1}}}.
\end{eqnarray*}
This finishes the proof.
\end{proof}

\begin{remark} \label{generalrealroots} 
Let $(h_0, \ldots, h_d) \in \RR^{d+1}$ such that
$h_i \geq 0$ for $0 \leq i \leq d$. An let $h_j' = \displaystyle{\sum_{i=0}^d} A(d+1,j,i+1) h_i$
for $0 \leq j \leq d$ then the polynomial $h_d' + h_{d-1}'t + \cdots h_0't^d$ has only
real and simple zeros.
\end{remark}
\begin{proof} 
The proof of Theorem \ref{realroots} never uses the fact that $h_i^\Delta$ is the
$h$-vector of a Boolean cell complex. The proof only relies on $h_i^\Delta \geq 0$
and the formula transforming $h_i^\Delta$ into $h_i^{\sd(\Delta)}$. Therefore,
the Remark follows from the proof of Theorem \ref{realroots}.
\end{proof}

Note that since $\sd(\Delta)$ is a flag simplicial complex; that is
all minimal non-faces are of size $2$; it is known that
$\hpol^{\sd(\Delta)} (t)$ has at least one real zero. This fact for flag 
simplicial complexes appears in \cite{ReiWel} or \cite{Gal}. As mentioned 
in \cite{ReiWel} this is actually a well known fact in commutative algebra, 
where it appears as a theorem saying that the numerator polynomial
of the Hilbert-series of a standard graded Koszul $k$-algebra has at least one
real zero. The following example shows that the assumptions of the theorem
are actually needed.   

\begin{example} \label{example-nonrealzeros}
Let $\Delta$ be the disjoint union of a $2$-simplex and a $1$-simplex.
Then  $\hpol^{\Delta}=(1,2,-3,1)$
and $\hpol ^{\sd (\Delta )}=(1,7,-3,1)$ so
$ \hpol^{\sd (\Delta )}(t)$ has only one real zero.
\end{example}

\begin{corollary} \label{CM-realrooted} 
Let $\Delta$ be a $(d-1)$-dimensional simplicial complex that is Cohen-Macaulay over some field $k$. Then
$$\hpol^{\sd(\Delta)}(t) = \sum_{i=0}^d h_i^{\sd(\Delta)} t^{d-i}$$ has only simple
and real zeros. In particular, 
$\hvec^{\sd(\Delta)} = (h_0^{\sd(\Delta)}, \ldots, h_d^{\sd(\Delta)})$ is a log-concave and
unimodal sequence.
In particular, all conclusions hold for Gorenstein, Gorenstein$^*$ simplicial complexes,
for simplicial
spheres and simplicial polytopes.
\end{corollary}
\begin{proof} 
All mentioned simplicial complex are Cohen-Macaulay over some field $k$. Therefore, it suffices to recall that
Cohen-Macaulay simplicial complexes have a non-negative $h$-vector without internal zeros.
This well known fact follows since if $\Delta$ is Cohen-Macaulay the sequence $(h_0^\Delta, \ldots,
h_d^\Delta)$ is the Hilbert-function $h_i^\Delta = \dim_k A_i$ of a $0$-dimensional standard graded 
$k$-algebra $A = A_0 \oplus \cdots \oplus A_d$ (see \cite[Chapter 5]{BH} for details). Being the dimension of a 
$k$-vectorspace clearly implies
that $h_i^\Delta \geq 0$. Now the fact that a standard graded $k$-algebra is generated in degree $1$
implies that there cannot be internal zeros in the sequence $(h_0^\Delta, \ldots,
h_d^\Delta)$.
\end{proof}

The preceding corollary assures that there is an abundance of interesting simplicial complexes 
satisfying the assumption of Theorem \ref{realroots}. At this point we would like to add
examples that show that the additional generality of Boolean cell complexes enriches the
class of complexes for which the theorem holds by interesting examples.

\begin{example} \label{injective-words}
An injective word over the alphabet $[n]$ is a word $i_1 \cdots i_j$ for which 
$i_r \neq i_s$ for all $1 \leq r < s \leq j$. We call $\{ i_1, \ldots, i_j \}$
the content of the injective word $i_1 \cdots l_j$. We order injective words
by saying $i_1' \cdots i_l' \preceq i_i \ldots i_j$ if there is a sequence
of indices $1 \leq s_1 < \cdots < s_l \leq j$ such that $i_1' \cdots i_l' =
i_{s_1} \cdots i_{s_l}$. For a simplicial complex $\Delta$ over ground set $[n]$
we denote by $\Gamma(\Delta)$ the partially ordered set of all injective words
with content in $\Delta$. It is easily seen that $\Gamma(\Delta)$ is the face poset
of a Boolean cell complex. In \cite{JoWe}, where the complexes $\Gamma(\Delta)$ 
are defined, it is shown that all $\Gamma(\Delta)$ have
non-negative $h$-vector. Thus providing a big class of Boolean cell complexes satisfying 
the assumptions of Theorem \ref{realroots}. It should be mentioned that if $\Delta$ is the
full simplex then the Boolean cell complex $\Gamma(\Delta)$ is a well studied object 
(see the references in \cite{JoWe}).  
\end{example}

At this point we are in position to relate our results to the Charney-Davis 
conjecture \cite{CharneyDavis}. The original conjecture states that if $\Delta$ is 
a flag $(d-1)$-dimensional Gorenstein$^*$ simplicial complex and $d$ is even then
$(-1)^{\lfloor \frac{d}{2} \rfloor} (h_0 -h_1 + \cdots +(-1)^d h_d) \geq 0$.  
Alternatively, the assertion can be phrased as 
$(-1)^{\lfloor \frac{d}{2}\rfloor} \hpol^{\sd(\Delta)}(-1) \geq 0$.
Since $d$ odd implies for Gorenstein$^*$ simplicial complexes by 
$h_i^\Delta = h_{d-i}^\Delta$ that $\hpol^{\sd(\Delta)} (-1) =0$
we will drop the assumption $d$ even and say that a simplicial complex $\Delta$
with reciprocal $h$-vector satisfies the Charney-Davis Conjecture if 
$(-1)^{\lfloor \frac{d}{2}\rfloor} \hpol^{\sd(\Delta)}(-1) \geq 0$.
The latter number is also called Charney-Davis quantity of $\Delta$.
Barycentric subdivisions already appear in one of the first instances where
the Charney-Davis conjecture was verified. Stanley \cite{st3} shows
that if $\Delta$ is the barycentric subdivision of the boundary complex of
a (not necessarily simplicial) polytope then the Charney-Davis conjecture holds
for $\Delta$. More recently, the work of Karu \cite{Ku} implies that for flag simplicial 
spheres that are order complexes of Gorenstein$^*$ posets the conjecture follows.
As a corollary of Theorem \ref{CM-realrooted} we obtain the 
following result which strengthens Karu's results in the case of 
subdivisions of Boolean cell complexes that are spheres. 

\begin{corollary} The Charney-Davis conjecture holds for
subdivisions of $(d-1)$-dimensional Boolean cell complexes $\Delta$ for which 
$h_i^\Delta \geq 0$ and $h_i^\Delta = h_{d-i}^\Delta$ for $0 \leq i \leq d$. 
If $d$ is even then the assumptions imply that the Charney-Davis quantity 
of $\sd(\Delta)$ is strictly positive.
\end{corollary}
\begin{proof} 
First let us verify the Charney-Davis conjecture for $\sd(\Delta)$ under the
given assumptions. 
By Corollary \ref{reciprocal} the assumptions imply that $\sd(\Delta)$ has 
a reciprocal $h$-polynomial $\hpol^{\sd(\Delta(t))}$. By Theorem \ref{realroots}
we also know that $\hpol^{\sd(\Delta)}$ has only real zeros. Since the coefficients of
$\hpol^{\sd(\Delta)}(t)$ are non-negative and $h_0^{\sd(\Delta)} = 1$ it follows 
that the zeros of $\hpol^{\sd(\Delta)}(t)$ are all strictly negative. Being a reciprocal
polynomial with non-zero constant coefficient also implies that if $\alpha$ is a zero 
then $\frac{1}{\alpha}$ is a zero. Thus the zeros are either $-1$ or come in 
pairs $\alpha < -1  < \frac{1}{\alpha} < 0$. If $-1$ is a zero then the the assertion
follows immediately. If $-1$ is not a zero then $d$ must be even and the $\frac{d}{2}$
factors $(-1 - \alpha) (-1 - \frac{-1}{\alpha})$ of $\hpol^{\sd(\Delta)}(-1)$ have all
negative sign. Thus $\hpol^{\sd(\Delta)}(-1)$ has sign $(-1)^{\frac{d}{2}}$, which implies
Charney-Davis. 

Let again $d$ be even. Since there is an even number of zeros of
$\hpol^{\sd(\Delta)}(t)$ which by Theorem \ref{realroots} are simple and real it follows
that $-1$ can only be a simple zero. But since all zeros come in pairs
$\alpha, \frac{1}{\alpha}$ this is impossible. Thus $\hpol^{\sd(\Delta)}(-1) \neq 0$ and
it follows from the first part of the proof that the Charney-Davis quantity is 
strictly positive. 

It is well know that Gorenstein$^*$ simplicial complexes fulfill the assumptions of the
corollary.
\end{proof}

Note that the arguments showing that a monic real reciprocal 
$h$-polynomial with only real zeros and non.negative coefficients satisfies the Charney-Davis 
conjecture can also be found in \cite{ReiWel}.

\begin{remark} 
It had been asked in \cite{CharneyDavis} 
whether for flag Gorenstein$^*$ simplicial complexes $\Delta$ the 
$h$-polynomial $\hpol^\Delta(t)$ has only real zeros.  This conjecture has been 
shown to fail by Gal \cite{Gal} for $\dim \Delta 
\geq 5$ even for the class of flag simplicial polytopes. For $d \leq 4$ an even more general 
assertion is true. Namely, in \cite[Corollary 4.14]{ReiWel} it is shown that if $A$ is a 
finitely generated Koszul Gorenstein $k$-algebra such that
the numerator $\hpol^A(t)$ 
of its Hilbert-series $\Hilb(A,t) = \frac{\hpol^A(t)}{(1-t)^d}$ has 
degree $\leq 3$ then $\hpol^A(t)$ has only real zeros. 
\end{remark}

\begin{question} Is the main result true for (not necessarily simplicial) 
polytopes and the  -- toric $h$-vector ? -- the usual $h$-vector ?
\end{question}

In the remaining section we study the effect of barycentric subdivision on 
the Stanley-Reisner ring of a simplicial complex. Recall, that for a 
simplicial complex $\Delta$ on ground set $\Omega$ the Stanley-Reisner ideal $I_\Delta$ is the 
ideal generated within $S = k[x_\omega~|~\omega \in \Omega]$ by 
$\xx_A = \prod_{\omega \in A} x_\omega$ for $A \not\in \Delta$. 
In particular, $I_\Delta$ is
minimally generated by $\xx_A$ for $A$ a minimal non-face of $\Delta$. 
By $k[\Delta]$ we denote the quotient $S/I_\Delta$. 
Since $I_\Delta$ is a homogeneous ideal, the quotient $k[\Delta]$ is a standard graded 
$k$-algebra $k[\Delta] \cong \bigoplus_{i \geq 0} A_i^\Delta$. 
By $\Hilb(k[\Delta],t) = \sum_{i \geq 0} \dim_k A_i^\Delta$ we denote the
Hilbert-series of $k[\Delta]$. 
In the process of passing from $\Delta$ to $\sd(\Delta)$ the ring
theoretic properties change as follows:
\begin{itemize}
\item[$\triangleright$] The Krull dimensions $\dim(k[\Delta])$ and $\dim(k[\sd(\Delta)])$
coincide.
\item[$\triangleright$] The Hilbert-series of $k[\Delta]$ and $k[\sd(\Delta)]$
are given as:
$$\Hilb(k[\Delta],t) = \frac{h_0^\Delta + \cdots + h_dt^t}{(1-t)^d}$$ and
$$\Hilb(k[\sd(\Delta)],t) = \frac{h_0^{\sd(\Delta)} + \cdots + h_d^{\sd(\Delta)}t^d}{(1-t)^d}.$$ 
\item[$\triangleright$] $k[\sd(\Delta)]$ is Koszul. 
\end{itemize}

\begin{question} Is there a version of the main result for graded rings ?
What is the ``barycentric subdivision'' of a standard graded $k$-algebra ?
For that we have in mind an operator $\sd(\cdot)$ that transforms a
standard graded algebra $A$ into a standard graded algebra $\sd(A)$ such
that:
\begin{itemize}
\item[$\bullet$] $\sd(k[\Delta]) = k[\sd(\Delta)]$ for simplicial complexes
$\Delta$. 
\item[$\bullet$] $\dim(A) = \dim(\sd(A))$.
\item[$\bullet$] The numerator polynomials of the Hilbert-series of $A$ and
$\sd(A)$ transform  ''similarly'' as the ones of $k[\Delta]$ and $k[\sd(\Delta)]$.
\item[$\bullet$] $\sd(A)$ is Koszul.
\end{itemize}
\end{question}

\section{Limiting Behavior}
\label{limiting-section}

For a number $n \geq 1$ and a Boolean cell complex $\Delta$ we denote by 
$\sd^n(\Delta)$ the result of an $n$-fold application of the subdivision 
operator $\sd$  to $\Delta$.
If we fix a geometric realization $|\Delta|$  of $\Delta$ in some real vectorspace then 
one may regard the geometric realizations $|\sd^n(\Delta)|$ of the complexes $\sd^n(\Delta)$ 
as triangulations of the same space with triangles of area converging to $0$. 

\begin{example}
If $\dim(\Delta) = 0$ then $sd(\Delta) = \Delta$ which 
clearly implies that the $h$-vector remains invariant under barycentric
subdivision. 
\end{example}

\begin{theorem} \label{limit}
For a number $d \geq 1$ there are real negative numbers 
$\alpha_1, \ldots, \alpha_{d-2}$ such that for every Boolean cell 
complex $\Delta$ of dimension $d-1$ there are sequences 
$(\beta_i^{(n)})_{n \geq 1}$, $1 \leq i \leq d$ of complex numbers such 
that

\begin{itemize}
\item[(i)] $\beta_i^n$, $1 \leq i \leq d$, are real for $n$ sufficiently large. 
\item[(ii)] $\displaystyle{\lim_{n \rightarrow \infty}} \beta_i^{(n)} = \alpha_i$ for $1 \leq i \leq d-3$.
\item[(iii)] $\displaystyle{\lim_{n \rightarrow \infty}} \beta_{d-1}^{(n)} = 0$.
\item[(iv)] $\displaystyle{\lim_{n \rightarrow \infty}} \beta_{d}^{(n)} = -\infty$.
\item[(v)] For $n \geq 1$ 
$$\prod_{i=1}^{d} (t-\beta_i^{(n)}) = \hpol^{\sd^n(\Delta)}(t).$$ 
\end{itemize}
\end{theorem}

Before we can come to a proof of Theorem \ref{limit} we need to analyze the
transformation from $\hvec^\Delta$ to $\hvec^{\sd(\Delta)}$ more
closely.

First we set up the matrices of the transformation that send $f$- and
$h$-vectors of simplicial complexes to those of their barycentric subdivision.
Let $\fmat_{d-1} = (f_{ij}^{d-1})_{0 \leq i,j \leq d}$ be the matrix with 
entries $f_{ij}^{d-1} := (i+1)!S(j,i+1)$.  
Let $\hmat_{d-1} = (h_{ij}^{d-1})_{0 \leq i,j \leq d}$ be the matrix with 
entries $h_{ij}^{d-1} := A(d+1,i,j+1)$. By Lemma \ref{f-subdivision}
and Theorem \ref{h-subdivision} can be rephrased as follows:

$$\fvec^{\sd(\Delta)} = \fmat_{d-1} \fvec^\Delta,$$ 

$$\hvec^{\sd(\Delta)} = \hmat_{d-1} \hvec^\Delta,$$ 

In the sequel we summarize some simple results on the matrices $\hmat_{d-1}$
and $\fmat_{d-1}$ that
will be used in the proof of the main theorem of this section.

\begin{lemma} 
Let $d \geq 1$.
Then:

\begin{itemize}
\item[(i)] The matrices $\fmat_{d-1}$ and $\hmat_{d-1}$ are similar.
\item[(ii)] The matrices $\fmat_{d-1}$ and $\hmat_{d-1}$ are diagonizable with
eigenvalue $1$ of multiplicity $2$ and eigenvalues $2!,3!, \ldots, d!$
or multiplicities $1$. In particular, there is a basis 
of rational eigenvectors.
\end{itemize}
\end{lemma}
\begin{proof} Since the transformation from $\fvec^\Delta$ to 
$\hvec^\Delta$ is an invertible linear transformation the matrices 
$\fmat_{d-1}$ and $\hmat_{d-1}$ are similar by 
Lemma \ref{f-subdivision} and Theorem \ref{h-subdivision}.
For the second assertion consider $\fmat_{d-1}$. Clearly,
$\fmat_{d-1}$ is an upper triangular matrix with diagonal 
$1,1,2!,3!, \ldots,d!$. Thus it remains to be shown that the 
eigenspace for eigenvalue $1$ is of dimension $2$. But this follows since the first and the
second unit vector are eigenvectors for the eigenvalue $1$. 
\end{proof}

The following lemmas give the crucial information on the eigenvectors of 
$\hmat_{d-1}$.

\begin{lemma} \label{feigenvecs} 
Let $d \geq 2$ and $v_1^1(d), v_1^2(d), v_2(d), \ldots, v_d(d)$ be a basis of
$\RR^d$ consisting of eigenvectors of the matrix $\fmat_{d-1}$, where 
$v_1^1(d), v_1^2(d)$ are eigenvectors
to the eigenvalue $1$ and $v_i(d)$ is an eigenvector for the eigenvalue $i!$, $2 \leq i \leq d$.
Then the vectors $v_1^1(d+1) := (v_1^1(d),0)$, $v_1^2(d+1) := (v_1^2(d),0)$ and
$v_i(d+1) := (v_i(d),0)$ are eigenvectors of $\fmat_d$ for the eigenvalues $1,1,2!, \ldots, d!$.
\end{lemma}
\begin{proof}
The assertion follows from the fact that $\fmat_{d-1}$ and $\fmat_d$ are upper triangular and
that if one deletes the $(d+1)$-st column and row from $\fmat_d$ then one obtains $\fmat_{d-1}$. 
\end{proof}

\begin{lemma} \label{hvecexpansion} Let $d \geq 2$ and 
let $w_1^1, w_1^2, w_2, w_3, \ldots, w_d$ be
a basis of eigenvectors of the matrix $\hmat_{d-1}$, where $w_1^1, w_1^2$ are eigenvectors 
to the eigenvalue $1$ and $w_i$ is the eigenvector for the eigenvalue $i!$, $2 \leq i \leq d$.
\begin{itemize}
\item[(i)] 
Let $\Delta$ be a $(d-1)$-dimensional simplicial complex.
If we expand $\hvec^\Delta = a_1^1 w_1^1 +a_1^2 w_1^2+ \sum_{i=2}^d a_iw_i$
in terms of the eigenvectors, then $a_d \neq 0$.
\item[(ii)] The first and the last coordinate entry in $w_2, \ldots, w_d$ is zero.
\item[(iii)] The vectors $w_1^1$ and $w_1^2$ can be chosen such that 
$$w_1^1 = (1,i_1, \ldots, i_{d-1},0) \mbox{~and~} w_1^2 = (0,j_1, \ldots, j_{d-1},1).$$
\item[(iv)] The vector $w_d$ can be chosen such that 
$w_d = (0,a_1, \ldots, a_{d-1},0)$ for
strictly positive rational numbers $a_i$, $1 \leq i \leq d-1$.
\end{itemize}
\end{lemma}
\begin{proof}
\begin{itemize}
\item[(i)]
{From} Lemma \ref{feigenvecs} we deduce that if one expands the f-vector of $\Delta$ in terms of
a basis of eigenvectors of $\fmat_{d-1}$ then by $f_{d-1}^\Delta \neq 0$ 
the coefficient of the eigenvector to the eigenvalue
$d!$ is non-zero. Since $\fmat_{d-1}$ and $\hmat_{d-1}$ are similar 
the assertion follows.  
\item[(ii)] From the definition of $\hmat_{d-1}$ is is clear that the first row is the
first and the last row the $(d+1)$-st unit vector. Thus if we denote by $\hmat_{d-1}'$ the
matrix which results from $\hmat_{d-1}$ by deleting the first and $(d+1)$-st rows and columns,
then the characteristic polynomial of $\hmat_{d-1}$ splits into $(1-t)^2$ times the 
characteristic polynomial of $\hmat_{d-1}'$. In particular, $\hmat_{d-1}'$ has
eigenvalues $2!, \ldots, d!$ and therefore is diagonizable. 
Let $w_2', \ldots, w_d'$ be eigenvectors of $\hmat_{d-1}'$ for the eigenvalues $2!, \ldots, d!$. Again by the fact
that the first and last row of $\hmat_{d-1}$ are the first and the $(d+1)$-st unit vector it follows
that $(0,w_2',0), \ldots, (0,w_d',0)$ are eigenvectors of $\hmat_{d-1}$ for the eigenvalues 
$2!, \ldots, d!$. Since all eigenspaces are of dimension $1$ the assertion follows.
\item[(iii)] Follows immediate from (ii).
\item[(iv)] Let $\hmat_{d-1}'$ be as in the proof of (ii). 
The entries column of $\hmat_{d-1}'$ are $A(d+1,j,i)$ for $1 \leq j \leq d-1$ and $2 \leq i \leq d$.
It is easily seen that these numbers are strictly positive.
Therefore, by the Perron-Frobenius Theorem \cite[Theorem 8.2.11]{horn} it follows that there 
is an eigenvector $w_d'$ for the eigenvector $d!$ with strictly positive entries.
Now the assertion follows by (ii).
\end{itemize}
\end{proof}

\begin{lemma} \label{topeigenvalue} Let $w_d = (0,a_1, \ldots, a_{d-1},0)$ be an eigenvector of
$\hmat_{d-1}$ for the eigenvalue $d!$ such that $a_i > 0$ for all $1 \leq i \leq d-1$.
Then the polynomial $a_{d-1}t + \cdots a_1t^{d-1}$ has only real and simple zero.
\end{lemma}
\begin{proof}
The vector $w_d$ satisfies the assumptions of Remark \ref{generalrealroots}. Therefore, for
$(h_0', \ldots, h_d') = \hmat_{d-1} w_d$ it follows that $h_d+h_{d-1}t+\cdots + h_0t^d$
has only real and simple zeros. Since $w_d$ is an eigenvector of $\hmat_{d-1}$ for the
eigenvalue $d!$ it follows that $h_0 = h_d = 0$ and $h_i = d! a_i$. Thus,
$d! (a_{d-1}t + \cdots a_1t^{d-1}) = h_d+h_{d-1}t+\cdots + h_0t^d$ and 
$a_{d-1}t + \cdots a_1t^{d-1}$ has only real and simple zeros.   
\end{proof}
 
We now formulate a technical lemma on the behavior of sequences of polynomials and its zeros
which will serve as the key ingredient in the proof of Theorem \ref{limit}. 

\begin{lemma} \label{technicallemma} 
Let $(g_n(t))_{n \geq 0}$ be a sequences of real polynomials of degree 
$d-2$, $f(t)$ another real polynomial of degree $d-2$ 
and $\rho > 1$ a real number such that:
\begin{itemize}
\item[$\triangleright$] $\lim_{n \rightarrow \infty} 1/\rho^n g_n(t) = 
0$, where the limit is taken in ${\bf R}^{d-1}$.
\item[$\triangleright$] Either all the coefficients of the polynomial $f(t)$ are strictly positive or
strictly negative and $f(t)$ has  
$d-2$ simple real zeros $- \infty < \gamma_1 < \cdots < \gamma_{d-2} < 0$.
\end{itemize}
Then for any real $h_d \neq 0$ there are real numbers $\alpha_i$, $1 \leq i \leq d-2$ and
sequences $(\beta_i^n)_{n \geq 0}$, $1 \leq i \leq d$ of complex numbers such that:
\begin{itemize}
\item[(i)] $\beta_i^n$, $1 \leq i \leq d$, are real for $n$ sufficiently large.
\item[(ii)] $\lim_{n \rightarrow \infty} \beta_i^n = \alpha_i$, $1 \leq i \leq d-2$.
\item[(iii)] $\lim_{n \rightarrow \infty} \beta_{d-1}^n = 0$.
\item[(iv)] $\lim_{n \rightarrow \infty} \beta_{d}^n = -\infty$.
\item[(v)] $\prod_{i = 0}^{d-1} (t-\beta_i^n) = 1 + tg_n(t) + tf(t) +t^d$. 
\end{itemize}
\end{lemma}
\begin{proof} 
Consider a zero $\gamma_i$ of the polynomial $f(t)$. Since $\gamma_i$ is a simple zero 
of $f_n(t)$ there is an $\epsilon > 0$ such that either 
(A) 
$f(x) > 0$ if $x \in [\gamma_i - \epsilon,\gamma_i)$
and $f(x) < 0$ if $x \in (\gamma_i,\gamma_i + \epsilon)$ or  
(B) 
$f(x) < 0$ if $x \in [\gamma_i - \epsilon,\gamma_i)$
and $f(x) > 0$ if $x \in (\gamma_i,\gamma_i + \epsilon)$.
Without loss of generality we may assume (A). Set 
$$R_n^- = \min_{x \in [\gamma_i - \epsilon,\gamma_i]} (1+xg_n(x) + \rho^nxf(x)+x^d)$$
$$R_n^+ = \max_{x \in [\gamma_i,\gamma_i+\epsilon]} (h_d + xg_n(x) + x\rho^nf(x) +x^d).$$
Since $\lim_{n \rightarrow \infty} 1/\rho^n g_n(t) = 0$ by assumption and $\lim_{n \rightarrow \infty} 1/\rho^n (1+t^d) = 0$ it 
follows that $\lim_{n \rightarrow \infty} 1/\rho^n R_n^- < 0$
and $\lim_{n \rightarrow \infty} 1/\rho^n R_n^- > 0$. Thus there is an $N(\epsilon)$ such that $R_n^- < 0$ and $R_n^+ > 0$ for $n \geq N(\epsilon)$.
In particular, there is a zero of $h_d + tg_n(t) + tf(t) +t^d$ in the interval $[\gamma_i - \epsilon,\gamma_i+\epsilon]$ for
$n > N(\epsilon)$. 
If we choose $\epsilon < 1/4 \min_{1 \leq i \leq d-3} (\gamma_{i+1} - \gamma_i)$ then the intervals
$[\gamma_i - \epsilon,\gamma_i+\epsilon]$ are disjoint and hence the zeros distinct. In this situation we 
denote by $\beta_i^n$ the zero of $1 + tg_n(t) + tf(t) +t^d$ in the interval $[\gamma_i - \epsilon,\gamma_i+\epsilon]$.
Our arguments show that at $\beta_i^n$ the polynomial $1 + tg_n(t) + tf(t) +t^d$ changes the sign in the same way as 
$tf(t)$ does at $\gamma_i$.  

We now have to distinguish between odd and even $d$ and strictly positive and strictly negative coefficients of $f(t)$.
Assume $d$ is even and $f(t)$ has only strictly positive coefficients.
{From} the assumption it follows that $xf(x) < 0$ for $x < \gamma_1$ and for 
$\gamma_{d-2} < x < 0$. Thus at $\beta_1^n$ the polynomial $1 + tg_n(t) + tf(t) +t^d$ changes sign from $-$ to $+$ and
at $\beta_{d-2}^n$ from $+$ to $-$. By $\lim_{x \rightarrow \pm \infty} (1 + tg_n(t) + tf(t) +t^d) = \infty$ it
follows that there is a zero $\beta_{d-1}^n > \beta_{d-2}^n$ and a zero $\beta_{d}^n < \beta_{1}^n$. 
If $d$ is odd or $f(t)$ has only strictly negative coefficients analogous arguments show the existence of zeros 
$\beta_{d-1}^n > \beta_{d-2}^n$ and $\beta_{d}^n < \beta_{1}^n$.

So far we have shown that for $n$ large enough all zeros of $h_d + tg_n(t) + tf(t) +t^d$ are simple and real. 

Since each sequence $(\beta_i^n)_{n > N(\epsilon}$,$1 \leq i \leq d-2$, lies in a compact interval and since we can choose 
the length of the interval to be arbitrarily small it follows that the sequences $(\beta_i^n)_{n > N(\epsilon}$,
$1 \leq i \leq d-2$, converge to fixed real numbers $\alpha_i$. Since the coefficients of $t^i$, $1 \leq i \leq d-1$, in the polynomial 
$1 + tg_n(t) + tf(t) +t^d$ converge to $\pm \infty$ by $\lim_{n \rightarrow \infty} 1/\rho^n g_n(t) = 0$, it follows that there must be 
a zero of large absolute value. On the other hand the fact that the constant term $h_d$ is independent of $n$ it follows that
there must also be a zero close to zero. Since $\beta_d^n$ is bounded from above by $\gamma_1$ it follows that 
$\lim_{n \rightarrow \infty} \beta_d^n = -\infty$. But this in turn implies $\lim_{n \rightarrow \infty} \beta_{d-1}^n = 0$.
\end{proof}

Now we are in position to provide a proof of Theorem \ref{limit}.

\begin{proof}[Proof of Theorem \ref{limit}] 
Let $w_1^1, w_1^2, w_2, \ldots, w_d$ be eigenvectors of $\hmat_{d-1}$ for the eigenvalues $1,1,2 \ldots, d!$
satisfying (iii) and (iv) of Lemma \ref{hvecexpansion}..
If we expand $\hvec^\Delta$ as a linear combination $\hvec^\Delta = a_1^1w_1^1+a_1^2w_1^2+a_2w_2 + \cdots a_dw_d$
then by Lemma \ref{hvecexpansion} we know that $a_d \neq 0$. The $h$-polynomial $\hpol^\Delta(t)$ decomposes
analogously into $a_1^1 H_1^1(t) + a_1^2H_1^2(t) + a_2 H_2(t) + \cdots + a_dH_d(t)$, where
$H_1^1(t), H_1^2(t), H_2(t), H_d(t)$ are the generating polynomials of the vectors $w_1^1, w_1^2, w_2, \ldots, w_d$
when read backwards. By Lemma \ref{hvecexpansion} (iii) it follows that 
$a_1^1 = h_d^\Delta$ and $a_1^2 = h_0^\Delta = 1$. Clearly, $\hpol^{\sd^n(\Delta}) = 
h_d^\Delta H_1^1(t) + H_1^2(t) + 2^{n-1}H_2 h_2(t) + \cdots + (d!)^na_dH_d(t)$. 
We set
$$g_n(t) := \frac{1}{t} (h_d^\Delta H_1^1(t) + H_1^2(t) + a_2 H_2(t) + \cdots + (d-1)!^{n-1}a_{d-1}H_{d-1}(t) - h_d^\Delta - t^d)$$
and
$$f(t) = 1/(d!t) \cdot a_dH_d(t).$$ 
By Lemma \ref{hvecexpansion} (ii)-(iv) is follows that $f(t)$ and $g_n(t)$ are polynomials of degree $d-2$. Moreover, 
either all coefficients of $f(t)$ are either strictly negative or strictly positive and 
$$\hpol^{\sd^n(\Delta)} = h_d + tg_n(t) +t(d!)^nf(t) +t^d.$$ 
By Lemma \ref{topeigenvalue} the polynomial $H_d(t)$ and therefore the polynomial $f(t)$ has only simple and real zeros 
that are all strictly negative. 
Thus $g_n(t)$ and $f(t)$ satisfy the assumption of Lemma \ref{technicallemma} and therefore the assertion follows
by Lemma \ref{technicallemma}.
\end{proof}

It seems to be a challenging problem to determine the numbers $\alpha_i$, $1 \leq i \leq d-2$, form Theorem \ref{limit}.
This in turn is closely related to determining eigenvectors of the matrices $\fmat_{d-1}$ and $\hmat_{d-1}$ for the eigenvalue 
$d!$. We were not able to give an explicit description of such an eigenvector.

\begin{problem} Give a description of an eigenvector of the matrices $\fmat_{d-1}$ and $\hmat_{d-1}$ for the eigenvalue
$d!$.
\end{problem}

\end{document}